\def\E{\end{document}}
\documentclass[10pt]{article}
\usepackage{amssymb,amsmath}

\begin{document}
\title{
 $H^s_x\times H^s_x$ scattering theory for a weighted gradient system of 3D radial defocusing NLS
}
  \author{Xianfa Song{\thanks{E-mail: songxianfa2004@163.com(X.F. Song)
 }}\\
\small Department of Mathematics, School of Mathematics, Tianjin University,\\
\small Tianjin, 300072, P. R. China
}

\maketitle
\date{}

\newtheorem{theorem}{Theorem}[section]
\newtheorem{definition}{Definition}[section]
\newtheorem{lemma}{Lemma}[section]
\newtheorem{proposition}{Proposition}[section]
\newtheorem{corollary}{Corollary}[section]
\newtheorem{remark}{Remark}[section]
\renewcommand{\theequation}{\thesection.\arabic{equation}}
\catcode`@=11 \@addtoreset{equation}{section} \catcode`@=12

\begin{abstract}

In this paper, using $I$-method, we establish $H^s_x\times H^s_x$ scattering theories for the following Cauchy problem
\begin{equation*}
\left\{
\begin{array}{lll}
iu_t+\Delta u=\lambda |v|^2u,\quad iv_t+\Delta v=\mu|u|^2v,\quad x\in \mathbb{R}^3,\ t>0,\\
u(x,0)=u_0(x),\quad v(x,0)=v_0(x),\quad x\in \mathbb{R}^3.
\end{array}\right.
\end{equation*}
Here $\lambda>0$, $\mu>0$, $(u_0,v_0)\in H^s_x(\mathbb{R}^3)\times H^s_x(\mathbb{R}^3)$ and $\frac{1}{2}<s<1$.



\end{abstract}

\section{Introduction}
\qquad In this paper, we consider the following Cauchy problem:
\begin{equation}
\label{02201}
\left\{
\begin{array}{lll}
iu_t+\Delta u=\lambda |v|^2u,\quad iv_t+\Delta v=\mu|u|^2v,\quad x\in \mathbb{R}^3,\ t>0,\\
u(x,0)=u_0(x),\quad v(x,0)=v_0(x),\quad x\in \mathbb{R}^3.
\end{array}\right.
\end{equation}
Here $\lambda>0$, $\mu>0$, $(u_0,v_0)\in H^s_x(\mathbb{R}^3)\times H^s_x(\mathbb{R}^3)$ and $\frac{1}{2}<s<1$. Model (\ref{02201}) often appears in in condensed matter theory, in quantum mechanics, in nonlinear optics, in plasma physics and in the theory of Heisenberg ferromagnet and magnons. $\lambda>0$ and $\mu>0$ mean the nonlinearities are defocusing ones. An interesting topic on (\ref{02201}) is scattering phenomenon which we will study in this paper. 

First, we would like to review some scattering results on the following Cauchy problem when $d=3$
\begin{equation}
\label{1104w1}
\left\{
\begin{array}{lll}
iu_t+\Delta u=\lambda|u|^pu,\quad x\in \mathbb{R}^d,\ t\in \mathbb{R},\\
u(x,0)=u_0(x),\quad x\in \mathbb{R}^d,
\end{array}\right.
\end{equation}
where $p>0$ and $\lambda\geq 0$, i.e., the equation in (\ref{1104w1}) has the defocusing nonlinearity. About the classical results on $L^2_x$, $\Sigma$ and $H^1_x$ scattering theories for (\ref{1104w1}), to see more details, we can refer to the books \cite{Bourgain19992, Cazenave2003, Dodson2019, Tao2006} and the numerous references therein. Recently, different types of scattering results on (\ref{1104w1}) when $d=3$ were established by many authors, we can see \cite{Beceanu??, Bourgain19991, Colliander2008, Dodson2012, Dodson20190, Duyckaerts2008, Gao2019, Kenig2010,  Killip2018, Killip2017, Killip2012, Killip2013, Killip2016, Killip20082, Masaki2019, Miao2013, Murphy2014, Murphy2015, Tao20071, Xie2013}. We specially point out that the $H^s_x$ scattering theory for (\ref{1104w1}) with $p=2$ is also studied by many authors. Bourgain established the $H^s_x$ scattering theory for (\ref{1104w1}) in \cite{Bourgain19981} with $p=2$ when $d=2$ and initial data in $H^s_x$ with $s>\frac{3}{5}$. In \cite{Bourgain19982}, for (\ref{1104w1}) with $p=2$ when $d=3$ and radial initial data $u_0$, Bourgain proved global well-posedness and scattering for $s>\frac{5}{7}$. In \cite{Colliander2002}, for (\ref{1104w1}) with $p=2$, through introducing an operator $I:H^s_x(\mathbb{R}^d)\rightarrow H^1_x(\mathbb{R}^d)$ and tracking  the change of $E(Iu(t))$, Colliander et. al. proved global well-posedness for the solution to (\ref{1104w1}) when $d=2$ for $s>\frac{4}{7}$ and when $d=3$ for $s>\frac{5}{6}$. Later, the result when $d=3$ was respectively extend to $s>\frac{4}{5}$ by Colliander et. al. in \cite{Colliander2004}, $s>\frac{5}{7}$ by Dodson in \cite{Dodson2013}, $s>\frac{2}{3}$ by Su in \cite{Su2012}. Recently, using $I$-method, Dodson obtained the
$H^s_x$ scattering theory for (\ref{1104w1}) with $p=2$ for $s>\frac{1}{2}$.

 There are some results on the scattering theory for a system of Schr\"{o}dinger equations. For $[H^1_x]^N$-scattering theory, we can refer to \cite{Cassano2015, Farah2017, Peng2019, Saanouni20191, Tarulli2019, Tarulli2016, Xu2016} and see the details. In \cite{Saanouni20192}, Saanouni established the $[H^s_x]^N$-scattering theory for $\frac{2}{3}<s<1$. Very recently, in \cite{Song1}, we studied
\begin{equation}
\label{826x1}
\left\{
\begin{array}{lll}
iu_t+\Delta u=\lambda |u|^{\alpha}|v|^{\beta+2}u,\quad iv_t+\Delta v=\mu|u|^{\alpha+2}|v|^{\beta}v,\quad x\in \mathbb{R}^d,\ t>0,\\
u(x,0)=u_0(x),\quad v(x,0)=v_0(x),\quad x\in \mathbb{R}^d.
\end{array}\right.
\end{equation}
Besides establishing the local well-posedness of the $H^1_x\times H^1_x$, $\Sigma\times \Sigma$ and $H^s_x\times H^s_x$ solutions, we found that there exists a critical exponents line $\alpha+\beta=2$ in the sense that (\ref{826x1}) always possesses a unique bounded $H^1_x\times H^1_x$-solution for any initial data $(u_0,v_0)\in H^1_x(\mathbb{R}^3)\times H^1_x(\mathbb{R}^3)$ if $\alpha+\beta\leq 2$, while there exist some initial data $(u_0,v_0)\in H^1_x(\mathbb{R}^3)\times H^1_x(\mathbb{R}^3)$ such that it doesn't have the unique global $H^1_x\times H^1_x$-solution if $\alpha+\beta>2$ and $\alpha=\beta$. We also established $H^1_x\times H^1_x$ and $\Sigma\times \Sigma$ scattering theories when $(\alpha,\beta)$ below the critical exponents line when $d=3$, and $\dot{H}^{s_c}_x\times \dot{H}^{s_c}_x$ scattering theory when $d\geq 5$, where $s_c=\frac{d-2}{2}>1$.

As a special case of (\ref{826x1}), (\ref{02201}) is a weighted gradient system of Sch\"{o}dinger equations. Therefore, we can define the following weighted
mass and energy
\begin{align}
M_w(u,v)=\mu\|u\|_{L^2_x}+\lambda\|v\|_{L^2_x},\quad E_w(u,v)=\frac{1}{2}\int_{\mathbb{R}^3}[\mu|\nabla u|^2+\lambda|\nabla v|^2+\lambda\mu|u|^2|v|^2]dx.\label{03222}
\end{align}

Naturally, we hope to establish the $H^s_x\times H^s_x$ scattering theory for (\ref{02201}) and $\frac{1}{2}<s<1$.
 
 In fact, using $I$-method, we will establish $H^s_x\times H^s_x$ scattering theory for (\ref{02201}) as follows.

{\bf Theorem 1.1($H^s_x\times H^s_x$ scattering theory for (\ref{02201})).} {\it
Assume that $\lambda>0$, $\mu>0$, $(u_0,v_0)$ radial and $(u_0,v_0)\in H^s_x(\mathbb{R}^3)\times H^s_x(\mathbb{R}^3)$, $\frac{1}{2}<s<1$. Then the initial value problem (\ref{02201}) is globally well-posedness and scattering, i.e., there exist $(u_+,v_+)\in H^s_x(\mathbb{R}^3)\times H^s_x(\mathbb{R}^3)$ and $(u_-,v_-)\in H^s_x(\mathbb{R}^3)\times H^s_x(\mathbb{R}^3)$ such that
\begin{align}
&\lim_{t\rightarrow +\infty}[\|u(t)-e^{it\Delta}u_+\|_{H^s_x(\mathbb{R}^3)}+\|v(t)-e^{it\Delta}v_+\|_{H^s_x(\mathbb{R}^3)}]=0,\label{02202}\\
&\lim_{t\rightarrow -\infty}[\|u(t)-e^{it\Delta}u_-\|_{H^s_x(\mathbb{R}^3)}+\|v(t)-e^{it\Delta}v_-\|_{H^s_x(\mathbb{R}^3)}]=0.\label{02203}
\end{align}
}

{\bf Remark 1.2:} 1. Scattering for (\ref{02201}) corresponds to $[\|u\|_{L^5_{t,x}(\mathbb{R}\times\mathbb{R}^3)}+\|v\|_{L^5_{t,x}(\mathbb{R}\times\mathbb{R}^3)}]<+\infty$.

 2. If $u_0(x)\equiv v_0(x)\equiv w_0(x)\in H^s_x(\mathbb{R}^3)$ and $\lambda=\mu$, then (\ref{02201}) degrades into the scalar Schr\"{o}dinger equation $iw_t+\Delta w=\lambda |w|^2w$ subject to $w(x,0)=w_0(x)$, and the conclusions of Theorem 1.3 meet with those of Theorem 1.2 in \cite{Dodson2019}. In this sense, our results cover those of Theorem 1.2 in \cite{Dodson2019}.

The rest of this paper is organized as follows.  In Section 2, we give some preliminaries. In Section 3, we  will use $I$-method to establish $H^s_x\times H^s_x$ scattering theory for (\ref{02201}).

\section{Preliminaries}
\qquad In this section, we will give some notations and useful lemmas which were mentioned in \cite{Dodson2019}.

We use $L^q_tL^r_x(I\times \mathbb{R}^3)$ to denote the Banach space of functions $u:I\times\mathbb{R}^3\rightarrow \mathbb{C}$ subject to the norm
$$
\|u\|_{L^q_tL^r_x(\mathbb{R}\times \mathbb{R}^3)}:=\left(\int_I\left(\int_{\mathbb{R}^3} |f(x)|^rdx\right)^{\frac{q}{r}}dt\right)^{\frac{1}{q}}<\infty
$$
for any spacetime slab $I\times \mathbb{R}^3$, with the usual modifications when $q$ or $r$ is infinity. Especially, we abbreviate $L^q_tL^r_x$ as $L^q_{t,x}$ if $q=r$.

The Fourier transform on $\mathbb{R}^3$ and the fractional differential operators $|\nabla |^s$ are defined by
$$
\hat{f}(\xi):=\int_{\mathbb{R}^3}e^{-2\pi ix\cdot \xi}f(x)dx,\quad \widehat{|\nabla|^sf}(\xi)=|\xi|^s\hat{f}(\xi),
$$
and the homogeneous and inhomogeneous Sobolev norms are respectively defined by
\begin{align}
\|f\|_{\dot{H}^s_x(\mathbb{R}^3)}:=\||\nabla|^sf\|_{L^2_x(\mathbb{R}^3)}=\||\xi|^s\hat{f}(\xi)\|_{L^2_x(\mathbb{R}^3)},\quad \|f\|_{H^s_x(\mathbb{R}^3)}:=\|(1+|\xi|^2)^{\frac{s}{2}}\hat{f}(\xi)\|_{L^2_x(\mathbb{R}^3)}.\label{03223}
\end{align}

Let $e^{it\Delta}$ be the free Schr\"{o}dinger propagator and the generated group of isometries $(\mathcal{J}(t))_{t\in \mathcal{R}}$. Then for $2\leq p\leq \infty$ and $\frac{1}{p}+\frac{1}{p'}=1$,
\begin{align}
&\mathcal{J}(t)f(x):=e^{it\Delta}f(x)=\frac{1}{(4\pi it)^{\frac{3}{2}}}\int_{\mathbb{R}^3}e^{i|x-y|^2/4t}f(y)dy,\quad t\neq 0,\label{829x1} \\
&\widehat{e^{it\Delta}f}(\xi)=e^{-4\pi^2it|\xi|^2}\hat{f}(\xi),\label{829x3}\\
& \|e^{it\Delta}f\|_{L^{\infty}_x(\mathbb{R}^3)}\lesssim |t|^{-\frac{3}{2}}\|f\|_{L^1_x(\mathbb{R}^3)},\quad   t\neq 0,\label{829x2}\\
&\|e^{it\Delta}f\|_{L^p_x(\mathbb{R}^3)}\lesssim |t|^{-3(\frac{1}{2}-\frac{1}{p})}\|f\|_{L^{p'}_x(\mathbb{R}^3)}, \quad   t\neq 0.\label{11101}
\end{align}
And the following Duhamel's formula holds
\begin{align}
u(t)=e^{i(t-t_0)\Delta}u(t_0)-i\int_{t_0}^te^{i(t-s)\Delta}(iu_t+\Delta u)(s)ds.\label{829x4}
\end{align}

We mention some definitions and estimates which also appeared in \cite{Dodson2019} below.

{\bf Definition 2.1(Strichartz space).} {\it Let $\tilde{S}^0(I)$ be the Strichartz space
\begin{align}
\tilde{S}^0(I)=L^{\infty}_tL^2_x(I\times\mathbb{R}^3)\cap L^2_tL^6_x(I\times\mathbb{R}^3),\label{022011}
\end{align}
and $\tilde{N}^0(I)$ be the dual
\begin{align}
\tilde{N}^0(I)=L^1_tL^2_x(I\times\mathbb{R}^3)+L^2_tL^{\frac{6}{5}}_x(I\times\mathbb{R}^3).\label{022012}
\end{align}
}

And the following Strichartz estimates hold
\begin{align}
\|e^{it\Delta}u_0\|_{S^0(\mathbb{R}\times \mathbb{R}^3)}\lesssim \|u_0\|_{L^2(\mathbb{R}^3)},\quad \|\int^t_0e^{i(t-\tau)\Delta}F(\tau)d\tau\|_{S^0(I\times\mathbb{R}^3)}\lesssim \|F\|_{N^0(I\times \mathbb{R}^3)}.\label{03225}
\end{align}

{\bf Definition 2.2(Littlewood-Paley decomposition).} {\it Given a radial and decreasing function $\psi\in C^{\infty}_0(\mathbb{R}^3)$ satisfying $\psi(x)=1$ for $|x|\leq 1$ and $\psi(x)=0$ for $|x|>2$. Let
\begin{align}
\phi_j(x)=\psi(2^{-j}x)-\psi(2^{-j+1}x)\label{02204}
\end{align}
and $P_j$ be the Fourier multiplier defined by
\begin{align}
\widehat{P_jf}(\xi)=\phi_j(\xi)\hat{f}(\xi)\quad {\rm for\ any}\ j.\label{02205}
\end{align}
This gives the following Littlewood-Paley decomposition
\begin{align}
f=\sum_{j=-\infty}^{\infty} P_jf \quad {\rm at\ least\ in \ the}\ L^2\ {\rm sense}.\label{02206}
\end{align}
}

In convenience, we write
\begin{align}
P_{\leq N} u=\sum_{j:2^j\leq N}P_ju,\quad P_{>N}u=u-P_{\leq N}u.\label{0220w1}
\end{align}
We may abbreviate $u_{\leq N}=P_{\leq N}u$ and $u_{>N}=P_{>N}u$ in the sequels.

The following proposition gives some well-known properties for the Littlewood-Paley decomposition.

{\bf Proposition 2.3.} {\it
For any $1<p<+\infty$,
\begin{align}
&\|f\|_{L^p_x(\mathbb{R}^3)}\thicksim_p \|\left(\sum_{j=-\infty}^{+\infty}|P_jf|^2\right)^{\frac{1}{2}}\|_{L^p_x(\mathbb{R}^3)},\label{02207}\\
&\|f\|_{L^p_x(\mathbb{R}^3)}\lesssim_s \|f\|_{\dot{H}^s_x(\mathbb{R}^3)}\quad {\rm for}\quad s=3(\frac{1}{2}-\frac{1}{p}).\label{02208}
\end{align}
We also have for $2\leq p\leq \infty$
\begin{align}
&\|P_jf\|_{L^p_x(\mathbb{R}^3)}\lesssim 2^{3j(\frac{1}{2}-\frac{1}{p})}\|P_jf\|_{L^2_x(\mathbb{R}^3)},\label{02209}
\end{align}
and the radial Sobolev embedding
\begin{align}
\||x|P_jf\|_{L^{\infty}(\mathbb{R}^3)}\lesssim \|P_jf\|_{\dot{H}^{\frac{1}{2}}_x(\mathbb{R}^3)}.\label{022010}
\end{align}
}

Let $\psi$ is the same $\psi$ as in Definition 1 and $\chi(x)=\psi(\frac{x}{2})-\psi(x)$. Then
\begin{align}
\|\psi(\frac{x}{R})e^{it\Delta}(P_ju_0)\|_{L^2_{t,x}(I\times \mathbb{R}^3)}\lesssim 2^{-\frac{j}{2}}R^{\frac{1}{2}}\|P_ju_0\|_{L^2_x(\mathbb{R}^3)},\label{022011}\\
\|\int_Ie^{-it\Delta}\psi(\frac{x}{R})(P_jF(t))dt\|_{L^2_x(\mathbb{R}^3)}\lesssim 2^{-\frac{j}{2}}R^{\frac{1}{2}}\|\psi(\frac{x}{R})P_jF\|_{L^2_{t,x}(\mathbb{R}^3)},\label{022012}
\end{align}
and for any $q<2$, if $F$ is supported on $|x|\leq R$,
\begin{align}
\||\nabla|^{1-\frac{1}{q}}\int_0^{\infty}e^{-it\Delta}F(t,x)dt\|_{L^2_x(\mathbb{R}^3)}\lesssim R^{1-\frac{1}{q}}\|F\|_{L^q_tL^2_x(\mathbb{R}\times\mathbb{R}^3)},\label{022013}\\
\lesssim R^{1-\frac{1}{q}}\|F\|_{X_R(I\times \mathbb{R}^3)},\label{022014}
\end{align}
where
\begin{align}
\|F\|_{X_R(I\times \mathbb{R}^3)}=R^{\frac{1}{q}-1}\|\psi(Rx)F\|_{L^q_tL^2_x(I\times \mathbb{R}^3)}+R^{\frac{1}{q}-1}\sum_{j\geq 0}2^{j(1-\frac{1}{q})}\|\chi(2^{-j}Rx)F\|_{L^q_tL^2_x(I\times \mathbb{R}^3)}.\label{022015}
\end{align}

{\bf Definition 2.4($U^p_{\Delta}$ spaces).} {\it
Let $1\leq p<\infty$ and $U^p_{\Delta}$ be an atomic space whose atoms are piecewise solutions of the linear equation,
\begin{align}
u_j=\sum_k1_{[t_k,t_{k+1})}e^{it\Delta}u_k,\quad \sum_k\|u_k\|^p_{L^2_x(\mathbb{R}^3)}=1.\label{0220x1}
\end{align}
Then for any $1\leq p<\infty$,
\begin{align}
\|u\|_{U^p_{\Delta}}=\inf\{\sum_j|c_j|:u=\sum_ju_j,\quad u_j \ {\rm are}\ U^p_{\Delta} \ atoms \}.\label{0220x2}
\end{align}
}

{\bf Proposition 2.5.} {\it If $u$ solves
\begin{align}
iu_t+\Delta u=F_1+F_2,\quad u(0,x)=u_0(x)\label{0220x3}
\end{align}
on the interval $0\in I\subset \mathbb{R}$, then for $q<2$,
\begin{align}
\||\nabla|^{1-\frac{1}{q}}u\|_{U^2_{\Delta}(I\times\mathbb{R}^3)}\lesssim_q \||\nabla|^{1-\frac{1}{q}}u_0\|_{L^2_x(\mathbb{R}^3)}+\|F_1\|_{X_R(I\times\mathbb{R}^3)}
+\||\nabla|^{1-\frac{1}{q}}F_2\|_{L^{2+}_tL^{6-}_x(I\times\mathbb{R}^3)}.\label{0220x4}
\end{align}
}

\section{Scattering result on (\ref{02201}) }
\qquad In this section, we will use the I-method to establish scattering result on (\ref{02201}).

{\bf Definition 3.1($I$-operator).} {\it Let $I:H^s_x(\mathbb{R}^3)\rightarrow H^1_x(\mathbb{R}^3)$
be the Fourier multiplier
\begin{align}
\widehat{If}(\xi)=m_N(\xi)\hat{f}(\xi),\label{0220x5}
\end{align}
where
\begin{equation}
\label{0220x6}
m_N(\xi)=\left\{
\begin{array}{lll}
1& {\rm if}\ |\xi|\leq N,\\
\frac{N^{1-s}}{|\xi|^{1-s}}& {\rm if}\ |\xi|\geq 2N.
\end{array}\right.
\end{equation}
}

We first give some results on the $I$-operator and the modified energy $E(Iu(t),Iv(t))$. By Sobolev embedding, we have
\begin{align}
&\quad E_w(Iu(t),Iv(t))=\frac{1}{2}\int_{\mathbb{R}^3}[\mu|\nabla Iu|^2+\lambda |\nabla Iv|^2+\lambda\mu|Iu|^2|Iv|^2]dx\nonumber\\
&\lesssim \int_{\mathbb{R}^3}[|\nabla Iu|^2+ |\nabla Iv|^2+|Iu|^4+|Iv|^4]dx\nonumber\\
&\lesssim \|Iu\|^2_{\dot{H}^1_x(\mathbb{R}^3)}+\|Iv\|^2_{\dot{H}^1_x(\mathbb{R}^3)}+\|Iu\|^2_{\dot{H}^1_x(\mathbb{R}^3)}\|u\|^2_{\dot{H}^{\frac{1}{2}}_x(\mathbb{R}^3)}
+\|Iv\|^2_{\dot{H}^1_x(\mathbb{R}^3)}\|v\|^2_{\dot{H}^{\frac{1}{2}}_x(\mathbb{R}^3)},\label{0220w3}
\end{align}
and consequently
\begin{align}
E_w(Iu(0),Iv(0))\lesssim C(\|u_0\|_{H^s_x},\|v_0\|_{H^s_x})N^{2(1-s)}.\label{0220w4}
\end{align}
Meanwhile,
\begin{align}
\|u(t)\|^2_{H^s_x(\mathbb{R}^3)}+\|v(t)\|^2_{H^s_x(\mathbb{R}^3)}\lesssim E_w(Iu(t),Iv(t))+M_w(Iu(t),Iv(t)).\label{0220w5}
\end{align}
Here
$$
M_w(Iu(t),Iv(t))=\mu\|u(t)\|^2_{L^2_x(\mathbb{R}^3)}+\lambda\|v(t)\|^2_{L^2_x(\mathbb{R}^3)}
$$

Moreover, we can use the rescaling
 \begin{align}
(u(t,x),v(t,x))\rightarrow (u_a(t,x),v_a(t,x))=a^{\frac{1}{2}}(u(a^2t,ax),v(a^2t,ax))\label{03091}
 \end{align}
 such that $E_w(Iu(t),Iv(t))\leq 1$. And there exists $a^{s-\frac{1}{2}}\thicksim C(\|u_0\|_{H^s_x(\mathbb{R}^3)}, \|v_0\|_{H^s_x(\mathbb{R}^3)})N^{s-1}$ such that
\begin{align}
&E_w(Iu_{\lambda}(0),Iv_{\lambda}(0))\leq \frac{1}{2},\label{0220w6}
\end{align}
and
\begin{align}
&\|u_{\lambda}(0)\|_{L^2_x(\mathbb{R}^3)}\lesssim C(\|u(0)\|_{H^s_x(\mathbb{R}^3)})N^{\frac{1-s}{2s-1}}\|u(0)\|_{L^2_x(\mathbb{R}^3)},\label{0220w7}\\
&\|v_{\lambda}(0)\|_{L^2_x(\mathbb{R}^3)}\lesssim C(\|v(0)\|_{H^s_x(\mathbb{R}^3)})N^{\frac{1-s}{2s-1}}\|v(0)\|_{L^2_x(\mathbb{R}^3)}.\label{0220w8}
\end{align}

{\bf Proposition 3.2(Weight coupled interaction Morawetz estimate).} {\it Assume that $(u,v)$ is a solution of (\ref{02201}) on some interval $J$. Then
\begin{align}
&\quad \||u|^4\|_{L^1_{t,x}(J\times \mathbb{R}^3)}+\||v|^4\|_{L^1_{t,x}(J\times \mathbb{R}^3)}+\||u|^2|v|^2\|_{L^1_{t,x}(J\times\mathbb{R}^3)}\nonumber\\
&\lesssim
[\|u\|^2_{L^{\infty}_tL^2_x(J\times\mathbb{R}^3)}
+\|v\|^2_{L^{\infty}_tL^2_x(J\times\mathbb{R}^3)}]
[\|u\|^2_{L^{\infty}_t\dot{H}^{\frac{1}{2}}_x(J\times\mathbb{R}^3)}
+\|v\|^2_{L^{\infty}_t\dot{H}^{\frac{1}{2}}_x(J\times\mathbb{R}^3)}].\label{0220w9}
\end{align}
}

{\bf Proof:} Similar to \cite{Song1}, we define the following weight-coupled interaction Morawetz potential:
\begin{align}
M^{\otimes_2}_a(t)&=2\mu^2\int_{\mathbb{R}^3}\int_{\mathbb{R}^3}\widetilde{\nabla} a(x,y)\Im[\bar{u}(t,x)\bar{u}(t,y)\widetilde{\nabla}(u(t,x)u(t,y))]dxdy\nonumber\\
&\quad+2\lambda^2\int_{\mathbb{R}^3}\int_{\mathbb{R}^3}\widetilde{\nabla} a(x,y)\Im[\bar{v}(t,x)\bar{v}(t,y)\widetilde{\nabla}(v(t,x)v(t,y))]dxdy\nonumber\\
&\quad+2\lambda\mu\int_{\mathbb{R}^3}\int_{\mathbb{R}^3}\widetilde{\nabla} a(x,y)\Im[\bar{u}(t,x)\bar{v}(t,y)\widetilde{\nabla}(u(t,x)v(t,y))]dxdy\nonumber\\
&\quad+2\lambda\mu\int_{\mathbb{R}^3}\int_{\mathbb{R}^3}\widetilde{\nabla} a(x,y)\Im[\bar{v}(t,x)\bar{u}(t,y)\widetilde{\nabla}(v(t,x)u(t,y))]dxdy,\label{829w1}
\end{align}
where $a(x,y)=|x-y|$, $\widetilde{\nabla}=(\nabla_x, \nabla_y)$, $x\in \mathbb{R}^3$ and $y\in \mathbb{R}^3$. By the inequality (4.10) in \cite{Song1}, we have
\begin{align*}
&\quad \int_I\int_{\mathbb{R}^3}[|u(t,x)|^4+|v(t,x)|^4]dxdt\lesssim |M^{\otimes_2}_a(t)|\nonumber\\
&\lesssim[\|u\|^2_{L^{\infty}_tL^2_x(J\times\mathbb{R}^3)}
+\|v\|^2_{L^{\infty}_tL^2_x(J\times\mathbb{R}^3)}]
[\|u\|^2_{L^{\infty}_t\dot{H}^{\frac{1}{2}}_x(J\times\mathbb{R}^3)}
+\|v\|^2_{L^{\infty}_t\dot{H}^{\frac{1}{2}}_x(J\times\mathbb{R}^3)}],
\end{align*}
which implies (\ref{0220w9}) because  $|u|^2|v|^2\leq |u|^4+|v|^4$. \hfill $\Box$

{\bf Lemma 3.3.}  {\it If $E(Iu(a_l),Iv(a_l))\leq 1$, $J_l=[a_l,b_l]$, and $\|u\|_{L^4_{t,x}(J_l\times \mathbb{R}^3)}+\|v\|_{L^4_{t,x}(J_l\times \mathbb{R}^3)}\leq \epsilon$ for some $\epsilon>0$ sufficiently small, then
\begin{align}
\|\nabla Iu\|_{\tilde{S}^0(J_l\times \mathbb{R}^3)}+\|\nabla Iv\|_{\tilde{S}^0(J_l\times \mathbb{R}^3)}\lesssim 1.\label{02211}
\end{align}
}

{\bf Proof:} Let
\begin{align}
Z_I(t)=\|\nabla Iu\|_{\tilde{S}^0([a_l,t]\times \mathbb{R}^3)}+\|\nabla Iv\|_{\tilde{S}^0([a_l,t]\times \mathbb{R}^3)}.\label{02212}
\end{align}

Applying $\nabla I$ to the equations of (\ref{02201}), we have
\begin{equation}
\label{02213}
\left\{
\begin{array}{lll}
i(\nabla Iu)_t+\Delta (\nabla Iu)=\lambda \nabla I(|v|^2u),\quad i(\nabla Iv)_t+\Delta (\nabla Iv)=\mu\nabla I(|u|^2v),\\
\nabla Iu(x,0)=\nabla Iu_0,\quad \nabla Iv(x,0)=\nabla Iv_0.
\end{array}\right.
\end{equation}
Using Strichartz estimates with $q'=r'=\frac{10}{7}$, then applying a fractional Leibniz rule, we obtain
\begin{align}
Z_I(t)&\lesssim \|\nabla Iu_0\|_{L^2_x(\mathbb{R}^3)}+\|\nabla Iv_0\|_{L^2_x(\mathbb{R}^3)}+\|\nabla I(|v|^2u)\|_{L^{\frac{10}{7}}_{t,x}([a_l,t]\times\mathbb{R}^3)}+\|\nabla I(|u|^2v)\|_{L^{\frac{10}{7}}_{t,x}([a_l,t]\times\mathbb{R}^3)}\nonumber\\
&\lesssim [\|\nabla Iu\|_{L^{\frac{10}{3}}_{t,x}([a_l,t]\times\mathbb{R}^3)}+\|\nabla Iv\|_{L^{\frac{10}{3}}_{t,x}([a_l,t]\times\mathbb{R}^3)}]
[\|u\|^2_{L^5_{t,x}([a_l,t]\times\mathbb{R}^3)}+\|v\|^2_{L^5_{t,x}([a_l,t]\times\mathbb{R}^3)}]\nonumber\\
&\quad+\|\nabla Iu_0\|_{L^2_x(\mathbb{R}^3)}+\|\nabla Iv_0\|_{L^2_x(\mathbb{R}^3)}.\label{02214}
\end{align}
Since the $L^{\frac{10}{3}}_{t,x}$ factors are bounded by $Z_I(t)$, we only need to consider the bounds of  $\|u\|^2_{L^5_{t,x}([a_l,t]\times\mathbb{R}^3)}$ and $\|v\|^2_{L^5_{t,x}([a_l,t]\times\mathbb{R}^3)}$. In fact, similar to the proof of (3.10) in \cite{Colliander2004}, we can get
\begin{align}
&\|u\|^2_{L^5_{t,x}([a_l,b_l]\times\mathbb{R}^3)}\lesssim \epsilon^{\delta_1}(\|\nabla Iu\|_{\tilde{S}^0([a_l,b_l]\times \mathbb{R}^3)})^{\delta_2},\label{02215'}\\
&\|v\|^2_{L^5_{t,x}([a_l,b_l]\times\mathbb{R}^3)}\lesssim \epsilon^{\delta_1}(\|\nabla Iv\|_{\tilde{S}^0([a_l,b_l]\times \mathbb{R}^3)})^{\delta_2} \label{02215}
\end{align}
for some $\delta_1>0$ and $\delta_2>0$.

By (\ref{02214}), (\ref{02215'}) and (\ref{02215}), we have
\begin{align}
Z_I(t)\lesssim 1+\epsilon^{\delta_3}(Z_I(t))^{1+\delta_4}\quad {\rm for \ some} \ \delta_3>0, \delta_4>0.\label{02216}
\end{align}
For $\epsilon$ sufficiently small, (\ref{02216}) yields (\ref{02211}).\hfill $\Box$

(\ref{02211}) implies that
\begin{align}
\|\nabla Iu\|_{U^2_{\Delta}(J_l\times \mathbb{R}^3)}+\|\nabla Iv\|_{U^2_{\Delta}(J_l\times \mathbb{R}^3)}\lesssim 1.\label{03092}
\end{align}

{\bf Proposition 3.4.} {\it Let $0\in J$ be an interval such that $E(Iu(t),Iv(t))\leq 1$ on $J$. Then for $N(s, \|u_0\|_{H^s_x}, \|v_0\|_{H^s_x})$ sufficiently large,
\begin{align}
\|P_{>\frac{N}{8}}\nabla Iu\|_{L^2_tL^6_x(J\times \mathbb{R}^3)}+\|P_{>\frac{N}{8}}\nabla Iv\|_{L^2_tL^6_x(J\times \mathbb{R}^3)}\lesssim 1.\label{02217}
\end{align}
}

{\bf Proof:} Decomposing
\begin{align}
|v|^2u&=\varnothing((v_{>\frac{M}{8}})^2u)+\varnothing(v_{>\frac{M}{8}}u_{>\frac{M}{8}}v_{\leq\frac{M}{8}})
+\varnothing(v_{\leq\frac{M}{8}}u_{\leq\frac{M}{8}}v_{>\frac{M}{8}})+\varnothing((v_{\leq\frac{M}{8}})^2u_{>\frac{M}{8}})\nonumber\\
&\quad+\varnothing((v_{\leq\frac{M}{8}})^2u_{\leq\frac{M}{8}}),\label{02218}\\
|u|^2v&=\varnothing((u_{>\frac{M}{8}})^2v)+\varnothing(u_{>\frac{M}{8}}v_{>\frac{M}{8}}u_{\leq\frac{M}{8}})
+\varnothing(u_{\leq\frac{M}{8}}v_{\leq\frac{M}{8}}u_{>\frac{M}{8}})+\varnothing((u_{\leq\frac{M}{8}})^2v_{>\frac{M}{8}})\nonumber\\
&\quad+\varnothing((u_{\leq\frac{M}{8}})^2v_{\leq\frac{M}{8}}),\label{02219}
\end{align}
we can get
\begin{align}
P_{>M}(|v|^2u)&=P_{>M}[\varnothing((v_{>\frac{M}{8}})^2u)]+P_{>M}[\varnothing(v_{>\frac{M}{8}}u_{>\frac{M}{8}}v_{\leq\frac{M}{8}})]
+P_{>M}[\varnothing(v_{\leq\frac{M}{8}}u_{\leq\frac{M}{8}}v_{>\frac{M}{8}})]\nonumber\\
&\quad+P_{>M}[\varnothing((v_{\leq\frac{M}{8}})^2u_{>\frac{M}{8}})],\label{022110}\\
P_{>M}(|u|^2v)&=P_{>M}[\varnothing((u_{>\frac{M}{8}})^2v)]+P_{>M}[\varnothing(u_{>\frac{M}{8}}v_{>\frac{M}{8}}u_{\leq\frac{M}{8}})]
+P_{>M}[\varnothing(u_{\leq\frac{M}{8}}v_{\leq\frac{M}{8}}u_{>\frac{M}{8}})]\nonumber\\
&\quad+P_{>M}[\varnothing((u_{\leq\frac{M}{8}})^2v_{>\frac{M}{8}})],\label{022111}
\end{align}
because $P_{>M}[\varnothing((v_{\leq\frac{M}{8}})^2u_{\leq\frac{M}{8}})]=0$ and $P_{>M}[\varnothing((u_{\leq\frac{M}{8}})^2v_{\leq\frac{M}{8}})]=0$.
Since $\nabla I$ is a Fourier multiplier whose symbol is increasing as $|\xi|\nearrow +\infty$, by the product rule and (\ref{0220x4}), we obtain
\begin{align}
&\quad \|\nabla IP_{>M}u(t)\|_{U^2_{\Delta}(J\times\mathbb{R}^3)}\nonumber\\
&\lesssim \|\nabla IP_{>M}u(0)\|_{L^2_x(\mathbb{R}^3)}+\|\nabla IP_{>M}\varnothing[(v_{>\frac{M}{8}})^2u]\|_{L^{2-}_tL^{\frac{6}{5}+}_x(J\times \mathbb{R}^3)}\nonumber\\
&\quad+\|\nabla IP_{>M}\varnothing[(v_{>\frac{M}{8}})(u_{>\frac{M}{8}})(v_{\leq\frac{M}{8}})]\|_{L^{2-}_tL^{\frac{6}{5}+}_x}+\|P_{>M}\varnothing[(v_{>\frac{M}{8}})(\nabla v_{\leq\frac{M}{8}})(u_{\leq\frac{M}{8}})]\|_{L^{2-}_tL^{\frac{6}{5}+}_x}\nonumber\\
&\quad+\|P_{>M}\varnothing[(v_{>\frac{M}{8}})(v_{\leq\frac{M}{8}})(\nabla u_{\leq\frac{M}{8}})]\|_{L^{2-}_tL^{\frac{6}{5}+}_x}
+\|P_{>M}\varnothing[(u_{>\frac{M}{8}})(\nabla v_{\leq \frac{M}{8}})(v_{\leq \frac{M}{8}})\|_{L^{2-}_tL^{\frac{6}{5}+}_x}\nonumber\\
&\quad+M^{\frac{1}{q}-1}\left(\|P_{>M}\varnothing[(\nabla Iv_{>\frac{M}{8}})(v_{\leq \frac{M}{8}})(u_{\leq \frac{M}{8}})]\|_{X_R}+\|P_{>M}\varnothing[(\nabla Iu_{>\frac{M}{8}})(v_{\leq \frac{M}{8}})^2]\|_{X_R}\right)\nonumber\\
&:=\|\nabla IP_{>M}u(0)\|_{L^2_x(\mathbb{R}^3)}+(I)+(II)+(III)+(IV)+(V)+(VI)+(VII)\label{02221}
\end{align}
if $M\leq N$. Here all the $L^{2-}_tL^{\frac{6}{5}+}_x$ norms are on $(J\times \mathbb{R}^3)$.

Choosing $\delta(\epsilon)>0$ such that $(2-\epsilon, \frac{6}{5}+\delta(\epsilon))$ is the dual of an admissible pair, using the properties of $\nabla I$ and H\"{o}lder's inequality, Berstein's inequality, we have
\begin{align}
(I)&\lesssim \|\nabla Iu\|_{L^{\infty-\epsilon}_tL^{2+\delta(\epsilon)}_x}\|P_{>\frac{M}{8}}v\|^2_{L^4_tL^6_x}+\|\nabla IP_{>\frac{M}{8}}v\|_{L^2_tL^6_x}
\|P_{>\frac{M}{8}}v\|_{L^{\infty}_tL^2_x}\|P_{\leq N} u\|_{L^{\infty-\epsilon}_tL^{6+\delta(\epsilon)}_x}\nonumber\\
&\quad+\|\nabla IP_{>\frac{M}{8}}v\|_{L^2_tL^6_x}\|P_{>\frac{M}{8}}v\|_{L^{\infty}_tL^3_x}\|P_{>N}u\|_{L^{\infty-\epsilon}_tL^{3+\delta(\epsilon)}_x}\nonumber\\
&\lesssim \|\nabla IP_{>\frac{M}{8}}v\|_{L^2_tL^6_x}\left(M^{-1}[\|\nabla Iu\|_{L^{\infty-\epsilon}_tL^{2+\delta(\epsilon)}_x}+\|P_{\leq N} u\|_{L^{\infty-\epsilon}_tL^{6+\delta(\epsilon)}_x}]+M^{-\frac{1}{2}}\|P_{> N} u\|_{L^{\infty-\epsilon}_tL^{3+\delta(\epsilon)}_x}\right)\nonumber\\
&\lesssim M^{-1}N^{\frac{3(1-s)}{2s-1}\cdot\frac{\epsilon}{2(2-\epsilon)}}\|\nabla I P_{>\frac{M}{8}}v\|_{U^2_{\Delta}(J\times \mathbb{R}^3)}.\label{02222}
\end{align}
Here all the norms are on $(J\times \mathbb{R}^3)$, and we have use the following facts
\begin{align}
&\|\nabla Iu\|_{L^{\infty-\epsilon}_tL^{2+\delta(\epsilon)}_x}+\|P_{\leq N} u\|_{L^{\infty-\epsilon}_tL^{6+\delta(\epsilon)}_x}\lesssim_{s, \|u_0\|_{H^s_x(\mathbb{R}^3)}, \|v_0\|_{H^s_x(\mathbb{R}^3)}} N^{\frac{3(1-s)}{2s-1}\cdot\frac{\epsilon}{2(2-\epsilon)}},\label{02223}\\
&\|P_{> N} u\|_{L^{\infty-\epsilon}_tL^{3+\delta(\epsilon)}_x}\lesssim_{s, \|u_0\|_{H^s_x(\mathbb{R}^3)}, \|v_0\|_{H^s_x(\mathbb{R}^3)}} N^{-\frac{1}{2}+\frac{3(1-s)}{2s-1}\cdot\frac{\epsilon}{2(2-\epsilon)}}.\label{02224}
\end{align}
Similarly,
\begin{align}
(II)\lesssim M^{-1}N^{\frac{3(1-s)}{2s-1}\cdot\frac{\epsilon}{2(2-\epsilon)}}[\|\nabla I P_{>\frac{M}{8}}u\|_{U^2_{\Delta}(J\times \mathbb{R}^3)}+\|\nabla I P_{>\frac{M}{8}}v\|_{U^2_{\Delta}(J\times \mathbb{R}^3)}]\label{02225}
\end{align}
and
\begin{align}
(III)&=\|P_{>M}\varnothing[(v_{>\frac{M}{8}})(\nabla v_{\leq\frac{M}{8}})(u_{\leq\frac{M}{8}})]\|_{L^{2-}_tL^{\frac{6}{5}+}_x}\lesssim \|\nabla v_{\leq \frac{M}{8}}\|_{L^{\infty-\epsilon}_tL^{2+\delta(\epsilon)}_x}\|v_{>\frac{M}{8}}\|_{L^2_tL^6_x}\|u_{\leq \frac{M}{8}}\|_{L^{\infty}_tL^6_x}\nonumber\\
&\quad\lesssim M^{-1}N^{\frac{3(1-s)}{2s-1}\cdot\frac{\epsilon}{2(2-\epsilon)}}[\|\nabla I P_{>\frac{M}{8}}u\|_{U^2_{\Delta}(J\times \mathbb{R}^3)}+\|\nabla I P_{>\frac{M}{8}}v\|_{U^2_{\Delta}(J\times \mathbb{R}^3)}],\label{02226}\displaybreak\\
(IV)&=\|P_{>M}\varnothing[(v_{>\frac{M}{8}})(v_{\leq\frac{M}{8}})(\nabla u_{\leq\frac{M}{8}})]\|_{L^{2-}_tL^{\frac{6}{5}+}_x}\lesssim \|\nabla u_{\leq \frac{M}{8}}\|_{L^{\infty-\epsilon}_tL^{2+\delta(\epsilon)}_x}\|v_{>\frac{M}{8}}\|_{L^2_tL^6_x}\|v_{\leq \frac{M}{8}}\|_{L^{\infty}_tL^6_x}\nonumber\\
&\quad\lesssim M^{-1}N^{\frac{3(1-s)}{2s-1}\cdot\frac{\epsilon}{2(2-\epsilon)}}[\|\nabla I P_{>\frac{M}{8}}u\|_{U^2_{\Delta}(J\times \mathbb{R}^3)}+\|\nabla I P_{>\frac{M}{8}}v\|_{U^2_{\Delta}(J\times \mathbb{R}^3)}],\label{02227}\\
(V)&=\|P_{>M}\varnothing[(u_{>\frac{M}{8}})(\nabla v_{\leq \frac{M}{8}})(v_{\leq \frac{M}{8}})]\|_{L^{2-}_tL^{\frac{6}{5}+}_x}\lesssim \|\nabla v_{\leq \frac{M}{8}}\|_{L^{\infty-\epsilon}_tL^{2+\delta(\epsilon)}_x}\|u_{>\frac{M}{8}}\|_{L^2_tL^6_x}\|v_{\leq \frac{M}{8}}\|_{L^{\infty}_tL^6_x}\nonumber\\
&\quad\lesssim M^{-1}N^{\frac{3(1-s)}{2s-1}\cdot\frac{\epsilon}{2(2-\epsilon)}}[\|\nabla I P_{>\frac{M}{8}}u\|_{U^2_{\Delta}(J\times \mathbb{R}^3)}+\|\nabla I P_{>\frac{M}{8}}v\|_{U^2_{\Delta}(J\times \mathbb{R}^3)}].\label{02228}
\end{align}

To give the estimates for (VI) and (VII), we first recall the following radial Sobolev embedding
\begin{align}
\||x|^{\frac{3}{2}-s}u\|_{L^{\infty}_x(\mathbb{R}^3)}\lesssim \|u\|_{\dot{H}^s_x(\mathbb{R}^3)},\quad
\||x|^{\frac{3}{2}-s}v\|_{L^{\infty}_x(\mathbb{R}^3)}\lesssim \|v\|_{\dot{H}^s_x(\mathbb{R}^3)}.\label{02229}
\end{align}
Interpolating (\ref{02229}), by Strichartz estimates, using (\ref{0220w7}), (\ref{0220w8}) and (\ref{03092}), we have
\begin{align}
&\|Iu\|^4_{L^4_tL^{\infty}_x(J\times \mathbb{R}^3)}+\|Iv\|^4_{L^4_tL^{\infty}_x(J\times \mathbb{R}^3)}\lesssim_{\|u_0\|_{H^s_x(\mathbb{R}^3)},\|v_0\|_{H^s_x(\mathbb{R}^3)}} N^{\frac{3(1-s)}{2s-1}},\label{0222x1}\\
&\||x|^{\frac{1}{2}}Iu\|_{L^{\infty-\epsilon}_tL^{\infty}_x}
+\||x|^{\frac{1}{2}}Iv\|_{L^{\infty-\epsilon}_tL^{\infty}_x}\lesssim_{s,\|u_0\|_{H^s_x(\mathbb{R}^3)},\|v_0\|_{H^s_x(\mathbb{R}^3)}} N^{\frac{(1-s)}{2s-1}\cdot[\frac{3\epsilon}{2(2-\epsilon)}+\frac{\epsilon}{(2-3\epsilon)}]}.\label{0222x2}
\end{align}
Choosing $R=N$ in (\ref{022011}), by (\ref{02229}) and (\ref{0222x2}), we get
\begin{align}
&\quad R^{\frac{1}{q}-1}M^{\frac{1}{q}-1}\|\psi(Rx)(\nabla IP_{>\frac{M}{8}}v)(v_{\leq \frac{M}{8}})(u_{\leq \frac{M}{8}})\|_{L^q_tL^2_x(J\times\mathbb{R}^3)}\nonumber\\
&\lesssim R^{\frac{1}{q}-1}M^{\frac{1}{q}-1}\|\psi(Rx)(\nabla IP_{>\frac{M}{8}}v)\|_{L^2_{t,x}(J\times\mathbb{R}^3)}\|v_{\leq\frac{M}{8}}\|^{\frac{2\epsilon}{(2-\epsilon)}}_{L^4_tL^{\infty}_x}
\|u_{\leq\frac{M}{8}}\|^{\frac{4-4\epsilon}{(2-\epsilon)}}_{L^{\infty}_{t,x}}\nonumber\\
&\lesssim_{s,\|u_0\|_{H^s_x(\mathbb{R}^3)},\|v_0\|_{H^s_x(\mathbb{R}^3)}} N^{\frac{-4+3\epsilon}{2(2-\epsilon)}}M^{\frac{-4+3\epsilon}{2(2-\epsilon)}}N^{\frac{3(1-s)}{2s-1}\cdot \frac{2\epsilon}{2-\epsilon}}M^{\frac{2-2\epsilon}{2-\epsilon}}\|\nabla IP_{>\frac{M}{8}}v\|_{U^2_{\Delta}(J\times \mathbb{R}^3)}.\label{0222x3}
\end{align}
Meanwhile, by (\ref{022011}) and (\ref{0222x2}),
\begin{align}
&\quad M^{\frac{1}{q}-1}\sum_{j\geq 0}R^{\frac{1}{q}-1}2^{j(1-\frac{1}{q})}\|\chi(2^{-j}Rx)(\nabla IP_{>\frac{M}{8}}v)(v_{\leq \frac{M}{8}})(u_{\leq \frac{M}{8}})\|_{L^q_tL^2_x(J\times\mathbb{R}^3)}\nonumber\\
&\lesssim M^{\frac{1}{q}-1}\sum_{j\geq 0}R2^{-j}R^{\frac{1}{q}-1}2^{j(1-\frac{1}{q})}\|\chi(2^{-j}Rx)(\nabla IP_{>\frac{M}{8}}v)\|_{L^2_{t,x}}\||x|^{\frac{1}{2}}Iv\|_{L^{\infty-\epsilon}_tL^{\infty}_x}\||x|^{\frac{1}{2}}Iu\|_{L^{\infty}_{t,x}}\nonumber\\
&\lesssim_{s,\|u_0\|_{H^s_x(\mathbb{R}^3)},\|v_0\|_{H^s_x(\mathbb{R}^3)}} N^{\frac{\epsilon}{2(2-\epsilon)}}M^{\frac{\epsilon}{2(2-\epsilon)}}N^{\frac{(1-s)}{2s-1}\cdot[\frac{3\epsilon}{2(2-\epsilon)}+\frac{\epsilon}{2-3\epsilon}]}
\|\nabla IP_{>\frac{M}{8}}v\|_{U^2_{\Delta}(J\times \mathbb{R}^3)}.\label{0222x4}
\end{align}
The estimate for (VI) can be given by (\ref{0222x3}) and (\ref{0222x4}). We can obtain the estimate for (VII) similar to (\ref{0222x3}) and (\ref{0222x4}).

Combining (\ref{02221}), (\ref{02222}), (\ref{02225}), (\ref{02226}), (\ref{02227}) with the estimates for (VI) and (VII), we can get
\begin{align}
\|\nabla IP_{>M}u\|_{U^2_{\Delta}(J\times\mathbb{R}^3)}&\lesssim_{s,\|u_0\|_{H^s_x(\mathbb{R}^3)},\|v_0\|_{H^s_x(\mathbb{R}^3)}, \epsilon} 1+\frac{N^{C_1(s)\epsilon}}{M^{1-C_2(s)\epsilon}}\|\nabla IP_{>\frac{M}{8}}u\|_{U^2_{\Delta}(J\times\mathbb{R}^3)}\nonumber\\
&\quad+\frac{N^{C_1(s)\epsilon}}{M^{1-C_2(s)\epsilon}}\|\nabla IP_{>\frac{M}{8}}v\|_{U^2_{\Delta}(J\times\mathbb{R}^3)}\label{0222x5}
\end{align}
for any $s$ with $\epsilon=\epsilon(s)$ sufficiently small, $C_1(s)\epsilon<\frac{1}{4}$ and $C_2(s)\epsilon<\frac{1}{4}$.

Similarly,
\begin{align}
\|\nabla IP_{>M}v\|_{U^2_{\Delta}(J\times\mathbb{R}^3)}&\lesssim_{s,\|u_0\|_{H^s_x(\mathbb{R}^3)},\|v_0\|_{H^s_x(\mathbb{R}^3)}, \epsilon} 1+\frac{N^{C_1(s)\epsilon}}{M^{1-C_2(s)\epsilon}}\|\nabla IP_{>\frac{M}{8}}u\|_{U^2_{\Delta}(J\times\mathbb{R}^3)}\nonumber\\
&\quad +\frac{N^{C_1(s)\epsilon}}{M^{1-C_2(s)\epsilon}}\|\nabla IP_{>\frac{M}{8}}v\|_{U^2_{\Delta}(J\times\mathbb{R}^3)}.\label{0222x6}
\end{align}

If $M>C(s,\|u_0\|_{H^s_x(\mathbb{R}^3)},\|v_0\|_{H^s_x(\mathbb{R}^3)})$, (\ref{0222x5}) and (\ref{0222x6}) mean that
\begin{align}
&\quad \|\nabla IP_{>M}u\|_{U^2_{\Delta}(J\times\mathbb{R}^3)}+\|\nabla IP_{>M}v\|_{U^2_{\Delta}(J\times\mathbb{R}^3)}\nonumber\\
&\lesssim_{s,\|u_0\|_{H^s_x(\mathbb{R}^3)},\|v_0\|_{H^s_x(\mathbb{R}^3)}} 1+N^{-\frac{1}{4}}C(s,\|u_0\|_{H^s_x(\mathbb{R}^3)},\|v_0\|_{H^s_x(\mathbb{R}^3)})^{-\frac{3}{4}}\nonumber\\
&\qquad\qquad \qquad\qquad \qquad\qquad\times [\|\nabla IP_{>\frac{M}{8}}u\|_{U^2_{\Delta}(J\times\mathbb{R}^3)}+\|\nabla IP_{>\frac{M}{8}}v\|_{U^2_{\Delta}(J\times\mathbb{R}^3)}].\label{0222w1}
\end{align}
Meanwhile, recalling (\ref{03092}), we have
\begin{align}
&\quad\|\nabla IP_{>C(s,\|u_0\|_{H^s_x},\|v_0\|_{H^s_x})N^{\frac{2}{3}}}u\|_{U^2_{\Delta}(J\times\mathbb{R}^3)}
+\|\nabla IP_{>C(s,\|u_0\|_{H^s_x},\|v_0\|_{H^s_x})N^{\frac{2}{3}}}v\|_{U^2_{\Delta}(J\times\mathbb{R}^3)}\nonumber\\
&\lesssim
N^{\frac{3(1-s)}{4s-2}}.\label{0222w2}
\end{align}
By induction, for $C(s,\|u_0\|_{H^s_x(\mathbb{R}^3)},\|v_0\|_{H^s_x(\mathbb{R}^3)})$ sufficiently large,
\begin{align}
\|\nabla IP_{>\frac{N}{8}}u\|_{U^2_{\Delta}(J\times \mathbb{R}^3)}+\|\nabla IP_{>\frac{N}{8}}u\|_{U^2_{\Delta}(J\times \mathbb{R}^3)}\lesssim_{s,\|u_0\|_{H^s_x},\|v_0\|_{H^s_x}} 1+N^{\frac{3(1-s)}{4s-2}}N^{-c\ln(N)}.\label{0222w3}
\end{align}
Choosing $N$ sufficiently large such that $c\ln(N)>\frac{3(1-s)}{4s-2}$, then using (\ref{0222w3}), we get
\begin{align}
\|\nabla IP_{>\frac{N}{8}}u\|_{U^2_{\Delta}(J\times \mathbb{R}^3)}+\|\nabla IP_{>\frac{N}{8}}u\|_{U^2_{\Delta}(J\times \mathbb{R}^3)}\lesssim_{s,\|u_0\|_{H^s_x},\|v_0\|_{H^s_x}} 1,\label{0222w4}
\end{align}
which implies (\ref{02217}). \hfill$\Box$

We give a bound on the modified energy increment.

{\bf Lemma 3.5.} {\it Let $N$ be sufficiently large such that
$$
\ln(N)\geq \frac{C_0(1-s)}{2s-1}+\ln(C(s,\|u_0\|_{H^s_x(\mathbb{R}^3)},\|v_0\|_{H^s_x(\mathbb{R}^3)})),
$$
then
\begin{align}
\int_J|\frac{d}{dt}E_w(Iu(t),Iv(t))|dt\lesssim \frac{1}{N^{1-}}.\label{02231}
\end{align}
}

{\bf Proof:} Note that
$$
E_w(Iu,Iv)=\frac{1}{2}\int_{\mathbb{R}^3}[\mu|\nabla Iu|^2+\lambda|\nabla Iv|^2-\lambda\mu|Iu|^2|Iv|^2]dx
$$
and $\Re((\overline{Iu_t})(iIu_t))=0$. Then
\begin{align}
&\quad \frac{d}{dt}E_w(Iu(t),Iv(t))\nonumber\\
&=\Re\int_{\mathbb{R}^3}\left\{\mu(\overline{Iu_t})[\lambda|Iv|^2(Iu)-\Delta Iu-iIu_t]
+\lambda(\overline{Iv_t})[\mu|Iu|^2(Iv)-\Delta Iv-iIv_t]\right\}dx\nonumber\\
&=\Re\int_{\mathbb{R}^3}\lambda\mu\left\{(\overline{Iu_t})[|Iv|^2(Iu)-I(|v|^2u)]+(\overline{Iv_t})[|Iu|^2(Iv)-I(|u|^2v)]\right\}dx.\label{02232}
\end{align}

Since $I$ is a Fourier multiplier which is constant in time and $\Delta$ commutes with $I$, by the equations of (\ref{02201}),
\begin{equation}
\label{02233}
\left\{
\begin{array}{lll}
iIu_t+\Delta Iu=\lambda |Iv|^2(Iu)+\lambda [I(|v|^2u)-|Iv|^2(Iu),\quad x\in \mathbb{R}^3,\ t>0,\\
iIv_t+\Delta Iv=\mu |Iu|^2(Iv)+\mu [I(|u|^2v)-|Iu|^2(Iv),\quad x\in \mathbb{R}^3,\ t>0.
\end{array}\right.
\end{equation}
Integrating the equations of (\ref{02233}) by parts, we have
\begin{align}
&\quad\frac{d}{dt}E_w(Iu(t),Iv(t))\nonumber\\
&=-\mu<i\nabla u, \nabla(|Iv|^2(Iu)-I(|v|^2u))>-\lambda\mu<iI(|v|^2u),(|Iv|^2(Iu)-I(|v|^2u))>\nonumber\\
&\quad -\lambda<i\nabla v, \nabla(|Iu|^2(Iv)-I(|u|^2v))>-\lambda\mu<iI(|u|^2v),(|Iu|^2(Iv)-I(|u|^2v))>.\label{02234}
\end{align}
Note that
\begin{align}
|Iv|^2(Iu)-I(|v|^2u)
&=(IP_{\leq\frac{N}{8}}v)^2(IP_{>\frac{N}{8}}u)-I((P_{\leq\frac{N}{8}}v)^2(P_{>\frac{N}{8}}u))\nonumber\\
&\quad+(IP_{>\frac{N}{8}}v)^2(IP_{\leq\frac{N}{8}}u)-I((P_{>\frac{N}{8}}v)^2(P_{\leq\frac{N}{8}}u))\nonumber\\
&\quad+(IP_{>\frac{N}{8}}v)^2(IP_{>\frac{N}{8}}u)-I((P_{>\frac{N}{8}}v)^2(P_{>\frac{N}{8}}u))\nonumber\\
&\quad+2(IP_{>\frac{N}{8}}v)(IP_{\leq\frac{N}{8}}v)(IP_{>\frac{N}{8}}u)-2I((P_{>\frac{N}{8}}v)(P_{\leq\frac{N}{8}}v)(P_{>\frac{N}{8}}u))\nonumber\\
&\quad+2(IP_{>\frac{N}{8}}v)(IP_{\leq\frac{N}{8}}v)(IP_{\leq\frac{N}{8}}u)-2I((P_{>\frac{N}{8}}v)(P_{\leq\frac{N}{8}}v)(P_{\leq\frac{N}{8}}u))\nonumber\\
&:=(1)+(2)+(3)+(4)+(5),\label{0223s1}\\
|Iu|^2(Iv)-I(|u|^2v)
&=(IP_{\leq\frac{N}{8}}u)^2(IP_{>\frac{N}{8}}v)-I((P_{\leq\frac{N}{8}}u)^2(P_{>\frac{N}{8}}v))\nonumber\\
&\quad+(IP_{>\frac{N}{8}}u)^2(IP_{\leq\frac{N}{8}}v)-I((P_{>\frac{N}{8}}u)^2(P_{\leq\frac{N}{8}}v))\nonumber\\
&\quad+(IP_{>\frac{N}{8}}u)^2(IP_{>\frac{N}{8}}v)-I((P_{>\frac{N}{8}}u)^2(P_{>\frac{N}{8}}v))\nonumber\\
&\quad+2(IP_{>\frac{N}{8}}u)(IP_{\leq\frac{N}{8}}u)(IP_{>\frac{N}{8}}v)-2I((P_{>\frac{N}{8}}u)(P_{\leq\frac{N}{8}}u)(P_{>\frac{N}{8}}v))\nonumber\\
&\quad+2(IP_{>\frac{N}{8}}u)(IP_{\leq\frac{N}{8}}u)(IP_{\leq\frac{N}{8}}v)-2I((P_{>\frac{N}{8}}u)(P_{\leq\frac{N}{8}}u)(P_{\leq\frac{N}{8}}v))\nonumber\\
&:=(6)+(7)+(8)+(9)+(10),\label{0223s1'}
\end{align}
and
\begin{align}
|m(\xi_2+\xi_3+\xi_4)-m(\xi_2)|\lesssim \frac{|\xi_3|+|\xi_4|}{|\xi_2|}.\label{0223s2}
\end{align}
By $E_w(Iu(t),Iv(t))\leq 1$ and (\ref{02217}), we can obtain
\begin{align}
&\quad \int_J<i\nabla Iu, \nabla((IP_{\leq\frac{N}{8}}v)^2(IP_{>\frac{N}{8}}u)-I((P_{\leq\frac{N}{8}}v)^2(P_{>\frac{N}{8}}u)))>dt\nonumber\\
&\quad+\int_J<i\nabla Iv, \nabla((IP_{\leq\frac{N}{8}}u)^2(IP_{>\frac{N}{8}}v)-I((P_{\leq\frac{N}{8}}u)^2(P_{>\frac{N}{8}}v)))>dt\nonumber\\
&\lesssim \frac{1}{N}[\|\nabla Iu\|_{L^{\infty}_tL^2_x}+\|\nabla Iv\|_{L^{\infty}_tL^2_x}][\|\nabla IP_{>\frac{N}{8}}u\|^2_{L^2_tL^6_x}+\|\nabla IP_{>\frac{N}{8}}v\|^2_{L^2_tL^6_x}]\nonumber\\
&\quad\times[\|Iu\|_{L^{\infty}_tL^6_x}+\|Iv\|_{L^{\infty}_tL^6_x}]\lesssim \frac{1}{N},\label{0223s3}\\
&\quad \int_J<i\nabla Iu, \nabla((IP_{>\frac{N}{8}}v)^2(IP_{\leq\frac{N}{8}}u)-I((P_{>\frac{N}{8}}v)^2(P_{\leq\frac{N}{8}}u)))>dt\nonumber\\
&\quad+\int_J<i\nabla Iv, \nabla((IP_{>\frac{N}{8}}u)^2(IP_{\leq\frac{N}{8}}v)-I((P_{>\frac{N}{8}}u)^2(P_{\leq\frac{N}{8}}v)))>dt\nonumber\\
&\lesssim [\|\nabla Iu\|_{L^{\infty}_tL^2_x}+\|\nabla Iv\|_{L^{\infty}_tL^2_x}][\|\nabla IP_{>\frac{N}{8}}u\|_{L^2_tL^6_x}+\|\nabla IP_{>\frac{N}{8}}v\|_{L^2_tL^6_x}]\nonumber\\
&\quad\times[\| IP_{>\frac{N}{8}}u\|_{L^2_tL^6_x}+\| IP_{>\frac{N}{8}}v\|_{L^2_tL^6_x}][\|P_{\leq \frac{N}{8}}u\|_{L^{\infty}_tL^6_x}+\|P_{\leq \frac{N}{8}}v\|_{L^{\infty}_tL^6_x}]\lesssim \frac{1}{N},\label{0223s4}\displaybreak\\
&\quad \int_J<i\nabla Iu, \nabla((IP_{>\frac{N}{8}}v)^2(IP_{>\frac{N}{8}}u)-I((P_{>\frac{N}{8}}v)^2(P_{>\frac{N}{8}}u)))>dt\nonumber\\
&\quad+\int_J<i\nabla Iv, \nabla((IP_{>\frac{N}{8}}u)^2(IP_{>\frac{N}{8}}v)-I((P_{>\frac{N}{8}}u)^2(P_{>\frac{N}{8}}v)))>dt\nonumber\\
&\lesssim [\|\nabla Iu\|_{L^{\infty}_tL^2_x}+\|\nabla Iv\|_{L^{\infty}_tL^2_x}][\|\nabla IP_{>\frac{N}{8}}u\|_{L^2_tL^6_x}+\|\nabla IP_{>\frac{N}{8}}v\|_{L^2_tL^6_x}]\nonumber\\
&\quad\times [\|P_{>\frac{N}{8}}u\|^2_{L^4_tL^6_x}]+\|P_{> \frac{N}{8}}v\|^2_{L^4_tL^6_x}]\lesssim \frac{1}{N},\label{0223s5}\\
&\quad \int_J < i\nabla Iu,\nabla[2(IP_{>\frac{N}{8}}v)(IP_{\leq\frac{N}{8}}v)(IP_{>\frac{N}{8}}u)-2I((P_{>\frac{N}{8}}v)(P_{\leq\frac{N}{8}}v)(P_{>\frac{N}{8}}u))]>dt\nonumber\\
&\quad+\int_J < i\nabla Iv,\nabla[2(IP_{>\frac{N}{8}}u)(IP_{\leq\frac{N}{8}}u)(IP_{>\frac{N}{8}}v)-2I((P_{>\frac{N}{8}}u)(P_{\leq\frac{N}{8}}u)(P_{>\frac{N}{8}}v))]>dt\nonumber\\
&\lesssim [\|\nabla Iu\|_{L^{\infty}_tL^2_x}+\|\nabla Iv\|_{L^{\infty}_tL^2_x}][\|\nabla IP_{>\frac{N}{8}}u\|_{L^2_tL^6_x}+\|\nabla IP_{>\frac{N}{8}}v\|_{L^2_tL^6_x}]\nonumber\\
&\quad\times[\| IP_{>\frac{N}{8}}u\|_{L^2_tL^6_x}+\| IP_{>\frac{N}{8}}v\|_{L^2_tL^6_x}][\|P_{\leq \frac{N}{8}}u\|_{L^{\infty}_tL^6_x}+\|P_{\leq \frac{N}{8}}v\|_{L^{\infty}_tL^6_x}]\lesssim \frac{1}{N},\label{0223x1}\\
&\quad \int_J < i\nabla Iu,\nabla[2(IP_{>\frac{N}{8}}v)(IP_{\leq\frac{N}{8}}v)(IP_{\leq\frac{N}{8}}u)-2I((P_{>\frac{N}{8}}v)(P_{\leq\frac{N}{8}}v)(P_{\leq\frac{N}{8}}u))]>dt\nonumber\\
&\quad+\int_J < i\nabla Iv,\nabla[2(IP_{>\frac{N}{8}}u)(IP_{\leq\frac{N}{8}}u)(IP_{\leq\frac{N}{8}}v)-2I((P_{>\frac{N}{8}}u)(P_{\leq\frac{N}{8}}u)(P_{\leq\frac{N}{8}}v))]>dt\nonumber\\
&\lesssim \frac{1}{N}[\|\nabla Iu\|_{L^{\infty}_tL^2_x}+\|\nabla Iv\|_{L^{\infty}_tL^2_x}][\|\nabla IP_{>\frac{N}{8}}u\|^2_{L^2_tL^6_x}+\|\nabla IP_{>\frac{N}{8}}v\|^2_{L^2_tL^6_x}]\nonumber\\
&\quad\times[\|Iu\|_{L^{\infty}_tL^6_x}+\|Iv\|_{L^{\infty}_tL^6_x}]\lesssim \frac{1}{N}.\label{0223x2}
\end{align}
All the norms above are on $J\times \mathbb{R}^3$.  (\ref{0223s3})--(\ref{0223x2}) give the estimates for $$-\mu<i\nabla u, \nabla(|Iv|^2(Iu)-I(|v|^2u))>-\lambda<i\nabla v, \nabla(|Iu|^2(Iv)-I(|u|^2v))>.$$

Now we consider the estimate for $$-\lambda\mu<iI(|v|^2u),(|Iv|^2(Iu)-I(|v|^2u))>-\lambda\mu<iI(|u|^2v),(|Iu|^2(Iv)-I(|u|^2v))>.$$

As a matter of convenience, we denote
\begin{align}
I(|v|^2)u&=(IP_{>\frac{N}{8}}v)^2(IP_{>\frac{N}{8}}u)+2(IP_{>\frac{N}{8}}v)(IP_{\leq \frac{N}{8}}v)(IP_{>\frac{N}{8}}u)+(IP_{\leq\frac{N}{8}}v)^2(IP_{>\frac{N}{8}}u)\nonumber\\
&\quad+(IP_{>\frac{N}{8}}v)^2(IP_{\leq \frac{N}{8}}u)+2(IP_{>\frac{N}{8}}v)(IP_{\leq \frac{N}{8}}v)(IP_{\leq\frac{N}{8}}u)+(IP_{\leq\frac{N}{8}}v)^2(IP_{\leq\frac{N}{8}}u)\nonumber\\
&:=(I)+(II)+(III)+(IV)+(V)+(VI),\label{02261}\\
I(|u|^2)v&=(IP_{>\frac{N}{8}}u)^2(IP_{>\frac{N}{8}}v)+2(IP_{>\frac{N}{8}}u)(IP_{\leq \frac{N}{8}}u)(IP_{>\frac{N}{8}}v)+(IP_{\leq\frac{N}{8}}u)^2(IP_{>\frac{N}{8}}v)\nonumber\\
&\quad+(IP_{>\frac{N}{8}}u)^2(IP_{\leq \frac{N}{8}}v)+2(IP_{>\frac{N}{8}}u)(IP_{\leq \frac{N}{8}}u)(IP_{\leq\frac{N}{8}}v)+(IP_{\leq\frac{N}{8}}u)^2(IP_{\leq\frac{N}{8}}v)\nonumber\\
&:=(I)'+(II)'+(III)'+(IV)'+(V)'+(VI)'.\label{02262}
\end{align}
Then
\begin{align}
&\quad \int_J <I(|v|^2u),(|Iv|^2(Iu)-I(|v|^2u))>dt\nonumber\\
&=\int_J\left(<(I),(1)>+<(I),(2)>+<(I),(3)>+<(I),(4)>+<(I),(5)>\right)dt\nonumber\\
&\quad+\int_J\left(<(II),(2)>+<(II),(3)>+<(II),(4)>+<(II),(5)>\right)dt\nonumber\displaybreak\\
&\quad+\int_J\left(<(III),(3)><(III),(4)>+<(III),(5)>+<(IV),(4)>\right)dt\nonumber\\
&\quad+\int_J\left(<(IV),(5)>+<(V),(5)>+<(VI),(1)>+<(VI),(2)>\right)dt\nonumber\\
&\quad+\int_J\left(<(VI),(3)>+<(VI),(4)>+<(VI),(5)>\right)dt.\label{02263}
\end{align}
Using Sobolev embedding theorem, Bernstein's inequality, we have
\begin{align}
&\quad\|(IP_{>\frac{N}{8}}v)^2(IP_{>\frac{N}{8}}u)\|_{L^2_{t,x}}\nonumber\\
&\lesssim [\|\nabla IP_{>\frac{N}{8}}u\|_{L^2_tL^6_x}+\|\nabla IP_{>\frac{N}{8}}v\|_{L^2_tL^6_x}]
[P_{>\frac{N}{8}}u\|^2_{L^{\infty}L^3_x}+P_{>\frac{N}{8}}v\|^2_{L^{\infty}L^3_x}]\lesssim \frac{1}{N},\label{0226x1}\\
&\int_J<(I),(1)>dt\lesssim \|(IP_{>\frac{N}{8}}v)^2(IP_{>\frac{N}{8}}u)\|^2_{L^2_{t,x}}\lesssim \frac{1}{N^2},\label{0226x2}
\end{align}
and
\begin{align}
\int_J<(I),(2)>dt&\lesssim \|(IP_{>\frac{N}{8}}v)^2(IP_{>\frac{N}{8}}u)\|_{L^2_{t,x}}[\|P_{>\frac{N}{8}}u\|_{L^2_tL^6_x}+\|P_{>\frac{N}{8}}v\|_{L^2_tL^6_x}]\nonumber\\
&\quad \times [\|P_{\leq \frac{N}{8}}u\|_{L^{\infty}_{t,x}}+\|P_{\leq \frac{N}{8}}v\|_{L^{\infty}_{t,x}}][\|P_{>\frac{N}{8}}u\|_{L^{\infty}_tL^3_x}+\|P_{>\frac{N}{8}}v\|_{L^{\infty}_tL^3_x}]\nonumber\\
&\lesssim \frac{1}{N^2},\label{0226x3}\\
\int_J<(I),(3)>dt&\lesssim \int_J[|P_{>\frac{N}{8}}u|^2+|P_{>\frac{N}{8}}v|^2][|P_{\leq\frac{N}{8}}u|^2+|P_{\leq\frac{N}{8}}v|^2][|u|^2+|v|^2]dt\nonumber\\
&\lesssim [\|P_{>\frac{N}{8}}u\|^2_{L^2_tL^6_x}+\|P_{>\frac{N}{8}}v\|^2_{L^2_tL^6_x}][\|P_{\leq \frac{N}{8}}u\|^4_{L^{\infty}_tL^6_x}+\|P_{\leq \frac{N}{8}}v\|^4_{L^{\infty}_tL^6_x}]\nonumber\\
&\quad+[\|P_{>\frac{N}{8}}u\|^4_{L^4_tL^6_x}+\|P_{>\frac{N}{8}}v\|^4_{L^4_tL^6_x}][\|P_{\leq \frac{N}{8}}u\|_{L^{\infty}_tL^6_x}+\|P_{\leq \frac{N}{8}}v\|^4_{L^{\infty}_tL^6_x}]\nonumber\\
&\lesssim \frac{1}{N^2},\label{0226x4}\\
\int_J<(VI),(3)>dt&\lesssim \frac{1}{N^2}[\|\nabla IP_{>\frac{N}{8}}u\|^2_{L^2_tL^6_x}+\|\nabla IP_{>\frac{N}{8}}v\|^2_{L^2_tL^6_x}][\|IP_{\leq \frac{N}{8}}u\|^3_{L^{\infty}_tL^6_x}+\|IP_{\leq \frac{N}{8}}v\|^3_{L^{\infty}_tL^6_x}]\nonumber\\
&\quad \times [\|IP_{\leq \frac{N}{8}}u\|_{L^{\infty}_tL^2_x}+\|IP_{\leq \frac{N}{8}}v\|_{L^{\infty}_tL^2_x}]\nonumber\\
&\lesssim \frac{1}{N^2}.\label{0227x1}
\end{align}
All the norms above are on $J\times \mathbb{R}^3$. Similar to (\ref{0226x3}), we can get $\int_J<(I),(4)>dt\lesssim \frac{1}{N^2}$. Similar to (\ref{0227x1}), we can obtain  $\int_J<(VI),(5)>dt\lesssim \frac{1}{N^2}$. Similar to (\ref{0226x4}), the estimates for other terms in (\ref{02263}) can be bounded by $\frac{1}{N^2}$ because they all contain two $P_{>\frac{N}{8}}$ factors and two $P_{\leq\frac{N}{8}}$ factors. Putting all the results above together, we obtain the bound for $\int_J <I(|v|^2u),(|Iv|^2(Iu)-I(|v|^2u))>dt$.

Similarly, the bound for $\int_J <I(|u|^2v),(|Iu|^2(Iv)-I(|u|^2v))>dt$ can be established. Then substituting all the results into (\ref{02234}), we get (\ref{02231}). \hfill $\Box$

Recalling that
\begin{align*}
\|u(t)\|^2_{L^2_x(\mathbb{R}^3)}+\|v(t)\|^2_{L^2_x(\mathbb{R}^3)}=\|u_0\|^2_{L^2_x(\mathbb{R}^3)}+\|v_0\|^2_{L^2_x(\mathbb{R}^3)},
\end{align*}
and
\begin{align*}
\|u(t)\|^2_{H^s_x(\mathbb{R}^3)}+\|v(t)\|^2_{H^s_x(\mathbb{R}^3)}\lesssim E_w(Iu(t),Iv(t))+\|u_0\|^2_{L^2_x(\mathbb{R}^3)}+\|v_0\|^2_{L^2_x(\mathbb{R}^3)},
\end{align*}
by Lemma 5.10 and (\ref{02231}), we have
\begin{align}
\|u(t)\|_{H^s_x(\mathbb{R}^3)}+\|v(t)\|_{H^s_x(\mathbb{R}^3)}\lesssim C(s,\|u_0\|_{H^s_x(\mathbb{R}^3)}, \|v_0\|_{H^s_x(\mathbb{R}^3)})[\|u_0\|_{H^s_x(\mathbb{R}^3)}+ \|v_0\|_{H^s_x(\mathbb{R}^3)}].\label{02281}
\end{align}

Let $(p,q)$ be a $\frac{1}{2}$-admissible pair which satisfies $\frac{2}{p}=3(\frac{1}{2}-\frac{1}{q}-\frac{1}{6})$. Interpolating (\ref{0220w9}) with (\ref{02281}), we can get a bounded on $L^p_tL^q_x$. Now we can partition $\mathbb{R}$ into finitely many pieces with $[\|u\|_{L^p_tL^q_x(J_l\times \mathbb{R}^3)}+\|v\|_{L^p_tL^q_x(J_l\times \mathbb{R}^3)}]<\epsilon$ and use a perturbation argument to obtain $[\|u\|_{L^5_{t,x}(\mathbb{R}\times\mathbb{R}^3)}
+\|v\|_{L^5_{t,x}(\mathbb{R}\times\mathbb{R}^3)}]<+\infty$, which implies scattering.\hfill $\Box$

\end{document}